\documentclass[11pt]{article}
\usepackage{tikz, amssymb, amsmath, graphicx, bm,indentfirst,parskip,titling,csquotes,float}
\usepackage{amsthm}
\usepackage[nottoc,notlot,notlof]{tocbibind}
\usepackage[twoside]{fancyhdr}
\usepackage[bookmarks=true,%
    colorlinks=true,%
    linkcolor=blue,
    citecolor=blue,%
    filecolor=blue,%
    menucolor=blue,%
    urlcolor=blue,%
    breaklinks=true]{hyperref}
\usetikzlibrary{braids}
\usepackage[english]{babel}
\usepackage{cite}
\usepackage{cleveref}
\crefname{lemma}{Lemma}{Lemmas}
\newtheorem{thm}{Theorem}[section]
\newtheorem{defn}[thm]{Definition}
\newtheorem{lem}[thm]{Lemma}

\newtheorem{prop}[thm]{Proposition}
\newtheorem{rmk}[thm]{Remark}

\pretitle{\centering\large\bf\MakeUppercase}
\title{The 1-Loop and Ohtsuki-Takata Invariants for Alternating Braids}
\posttitle{\par\vskip 0.5em}

\author{Alexander R. Pierson \thanks{piersoal@oregonstate.edu}}

\date{}

\begin{document}

\maketitle

\begin{abstract}
    In a step towards resolving the 1-loop conjecture, Ohtsuki and Takata defined their own eponymous invariant in terms of an open knot diagram, and showed that it is equal to the adjoint Reidemeister torsion for two-bridge knots. Garoufalidis and Yoon then showed that the 1-loop invariant is equal to the Ohtsuki-Takata invariant using a technique involving a collapsed octahedral decomposition. This paper expands those techniques to relate the matrix determinants of the 1-loop and Ohtsuki-Takata invariants for the case of reduced alternating braids.
\end{abstract}

\tableofcontents
\section{Introduction}

Kashaev's volume conjecture \cite{kashaev1995link} is one of the great bridges linking 3-manifold topology and knot theory-the conjecture that Kashaev's eponymous invariant, shown to be equal to an evaluation of the colored Jones polynomial \cite{murakami2001colored}, has a growth rate related to the volume of the knot complement. 
$$\lim_{N\to\infty} \frac{2\pi}{N} log|\langle K\rangle _N| =vol(S^3\setminus K).$$
Since then, this conjecture has been generalized, and Dimofte and Garoufalidis \cite{dimofte2013quantum} introduced a conjectural power series for the Kashaev invarant associated to an essential ideal triangulation of a cusped hyperbolic 3-manifold, and proved that its constant term is a topological invariant. In the same paper, they conjectured that this so called \textit{1-loop invariant} is equal to the adjoint Reidemeister torsion.
\par
In \cite{ohtsuki2015kashaev}, Ohtsuki and Takata defined their own eponymous invariant, and showed that it is equal to the adjoint Reidemeister torsion in the case of 2-bridge knots. Garoufalidis and Yoon\cite{garoufalidis20261} used a calculation coming from a collapsed octahedral decomposition to show show equivalence in the 2-bridge case of the Ohtsuki-Takata and 1-loop invariants, resolving the 1-loop conjecture for 2-bridge knots. 
\par
The 1-loop conjecture has been resolved, proved additionally for fibered 3-manifolds \cite{dunfield20231}, fundamental shadow links \cite{pandey2023geometry}, and finally in general\cite{ming2026adjoint}, but the problem of whether the Ohtsuki-Takata invariant is equal to the 1-loop invariant in general is unresolved.
\par
In this paper we will expand the techniques of Garoufalidis and Yoon to the case of reduced alternating braid diagrams of hyperbolic knots, showing a relation between the determinants of the matrices in each invariant. 
\par
It is worth noting that not all alternating hyperbolic knots have alternating braid diagrams. For example, the $5_2$ knot is alternating and hyperbolic, but does not have an alternating braid diagram, as shown in \cite{guevara-hernandez2021kawauchi}. Nonetheless, this is a large class of knots, and includes some knots that are not two-bridge knots, such as $8_{16}$. Additionally, the restriction to this kind of diagram is beneficial. The 1-loop invariant is defined in terms of an essential triangulation, which is to say that no edge loop is null-homotopic\cite{hodgson2015triangulations}. This was shown to be equivalent to the condition it has a solution to the gluing equations whose holonomy representation is discrete and faithful (a \textit{geometric} solution)\cite{segerman2011pseudo,tillmann2012degenerations}. Alternating reduced diagrams of alternating hyperbolic knots induce a corresponding ideal triangulation which is essential \cite{garoufalidis2016non,sakuma2018application}, so we are able to define the 1-loop invariant.
\par
The structure of this paper is as follows. Section 2 is dedicated to preliminary results and definitions, including a description of a triangulation obtained from an octahedral decomposition, as well as the formal definition of the 1-loop and Ohtsuki-Takata invariants. Section 3 is work on the case of reduced alternating braids, for which we discuss the particular method of labeling variables of the Ohtsuki-Takata invariant and eliminating shape parameters in the 1-loop invariant, and a computation of the corresponding determinants.

\subsection*{Acknowledgments} The author would like to thank his advisor, Ren Guo, for his continued support and insight, and Seokbeom Yoon for helpful discussions.

\section{Preliminaries}
\subsection{Collapsed Octahedral Decomposition}
A primary idea of relevance to this paper is the construction of an ideal triangulation of a knot complement which comes from collapsing an octahedral decomposition along an open knot diagram. This is a classical and well-known result due to Thurston, and described in detail in \cite{weeks2005computation}. We will describe this here, but see also \cite{garoufalidis20261} for more details of the construction. 
\par 
Given a knot $K$ in $S^3$ and a planar diagram $D$ of $K$, we define a \textit{segment} to be the part of the knot in the diagram that lies between two crossings, an \textit{over-arc} to be a series of consecutive segments which meet exactly two crossings in an under-pass, and an \textit{under-arc} to be a series of consecutive segments which meet exactly two crossings in an over-pass. In the alternating case, over/under-arcs thus consist of two consecutive segments. 
\par
\begin{figure}[H]
	\centering
	\includegraphics[width=.3\textwidth]{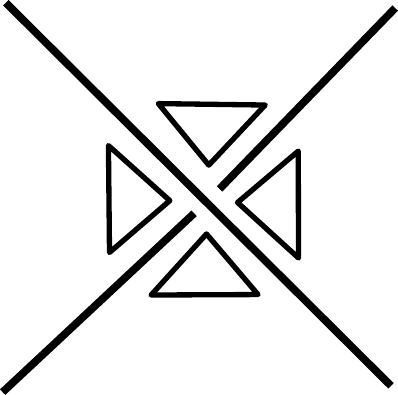}
	\caption{An Octahedron at a Crossing}
	\label{fig:CrossOct}
\end{figure}

 
Place an ideal octahedron at each crossing, where the octahedron can be decomposed into four glued tetrahedra, each of which corresponds to a corner in $D$. Seen from above this looks as in Figure \ref{fig:CrossOct}.
\par 
The edges of each octahedron can be classified into four types:
\begin{itemize}
	\item C-edges, which join the top and bottom vertices, and which correspond to a crossing of $D$,
	\item R-edges, which do not meet either the top or bottom vertex, and which correspond to a region of $D$,
	\item O-edges, which are transverse to the over-pass at the crossing, when viewed from the top, and correspond to an over-arc(not necessarily the over-arc containing the over-pass of the crossing),
	\item U-edges, which are transverse to the under-pass at the crossing, when viewed from the top, and correspond to an under-arc(not necessarily the under-arc containing the under-pass of the crossing).
\end{itemize}
Faces of octahedra are paired along segments, with two faces of the octahedron at one end of the segment being paired with two faces of the octahedron at the other end. The faces are the two upper faces if the segment meets the crossing at an over-pass, and the two lower faces if the segment meets the crossing at an under-pass. Face-pairing is done so that edges of the same type are identified. In the alternating case, upper faces are always paired with lower faces. 
\par 
This decomposition is an ideal triangulation of the knot complement minus two points, which we call $\pm \infty$. The upper and lower vertices of each octahedron meet the ideal vertex of the knot, and opposite equatorial vertices of each octahedron will meet the vertices at $\pm \infty$. For more details of this octahedral decomposition, as well as some illustrative figures, see \cite{kim2018octahedral}. 
\par
\begin{figure}[H]
	\centering
	\includegraphics[width=.75\linewidth]{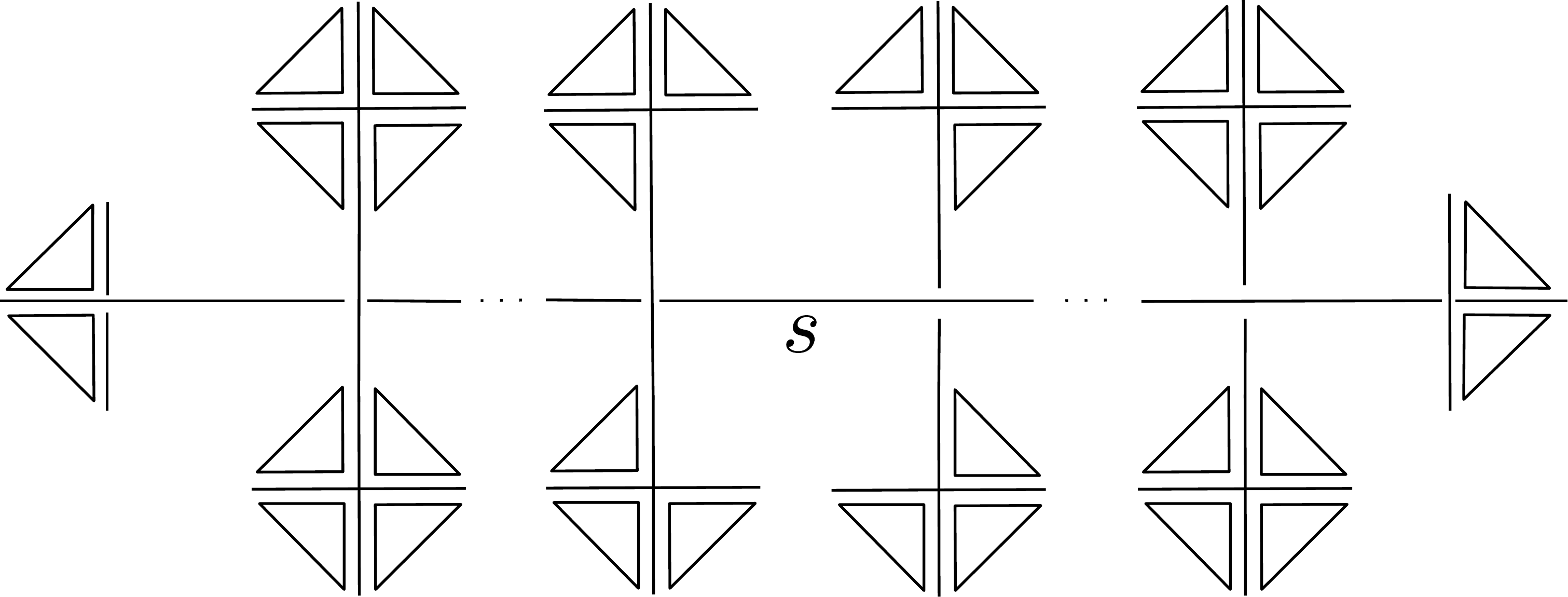}
	\caption{Collapsing the Octahedral Decomposition}
	\label{fig:Deletions}
\end{figure}

To turn this decomposition into an ideal triangulation of the knot complement, fix an alternating segment of $D$ (in the alternating case then we can choose any segment), and call is $s$. Then we remove all tetrahedra along the over-pass $o$ and under-pass $u$ containing $s$, and at corners in the regions adjacent to $s$. See Figure \ref{fig:Deletions}. 
Then interior faces are paired for octahedra along the regions adjacent to $s$, and where an exterior face would have been paired with the face of a tetrahedron that was removed, it is now paired with the face of a tetrahedron which also had a face paired to one of the removed tetrahedron. More specifically, call the crossing at the endpoint of $u$ which is not an endpoint of $s$ $c_1$, and call the crossing at the endpoint of $o$ which is not an endpoint of $s$ $c_m$, and enumerate all crossings in between $c_2,...,c_{m-1}$. For $i\in \{2,...,m-1\}$, there are two octahedra that would have paired faces with faces of the octahedron at $c_i$, which pair with each other instead. Two upper/lower faces would have been paired with two lower faces of the octahedron at $c_1$, but are instead paired with interior faces of the tetrahedra that remain at $c_1$. Similarly, two upper/lower faces would have been paired with two upper faces of the octahedron at $c_m$, but are instead paired with interior faces of the tetrahedra that remain at $c_m$.
\par Edges of different type can be identified under this new face-pairing, and in this case their type changes such that any O/U-edges retain their type (and thus can be both an O-edge and a U-edge), and an R-edge identified with a C-edge becomes an R-edge. See \cite{garoufalidis20261} for a precise description of edge identifications. 
\par Fixing the segment $s$ of $D$ allows us to turn $D$ into $\mathring{D}$, the open diagram obtained from $D$ by cutting along $s$. We can see that after collapsing, the R-edges corresponding to the unbounded regions as well as any regions adjacent to $o$ or $u$ disappear, as do any C-edges corresponding to crossings which meet $o$ or $u$. All other R and C-edges remain. 
\par
We will say that an open diagram is \textit{essential} if the corresponding collapsed triangulation is essential, which is to say the gluing and completeness equations have a geometric solution, one whose holonomy representation is discrete and faithful.

\subsection{The 1-loop Invariant}
Let \textit{M} be an oriented hyperbolic 3-manifold with one cusp, and let $\mathcal{T}$ be an oriented regular ideal triangulation of \textit{M}. By \textit{regular}, we mean a triangulation that recovers the complete hyperbolic structure. For a more complete discussion of regularity, see \cite{dimofte2013quantum}.
\par 
To each ideal tetrahedron we associate a set of \textit{shape parameters}, assigned to each edge such that opposite edged of tetrahedra have the same shape parameter, and $$z'=\frac{1}{1-z}\quad z''=1-\frac{1}{z}.$$
\begin{figure}[H]
	\centering
	\includegraphics[width=.4\linewidth]{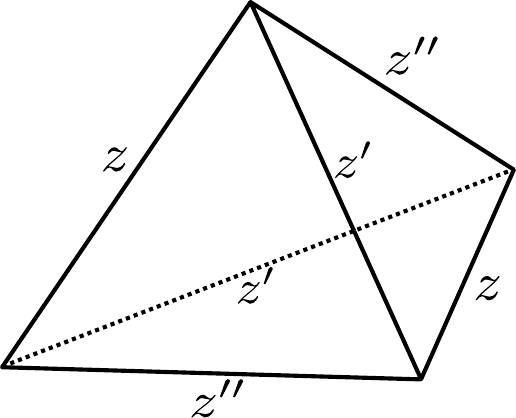}
	\caption{Shape Parameters of a Tetrahedron}
	\label{fig:spams}
\end{figure}
\begin{defn}[Dimofte-Garoufalidis \cite{dimofte2013quantum}]
Given \textit{M} and $\mathcal{T}$, an \textit{enhanced Neumann-Zagier datum} \text is the quintuple $$\hat{\beta}_{\mathcal{T}}=(\mathnormal{z},G,G',G'',\bm{f})$$ where:
\begin{itemize}
    \item $\mathnormal{z}=(\mathnormal{z_1,...,z_N)}$ is a solution to the gluing equations.
    \item $G$, $G'$ and $G''$ are the \textit{Neumann-Zagier matrices}, matrices with integer entries that encode the gluing information.
    \item $\bm{f}$ is a combinatorial flattening of $\mathcal{T}$.
\end{itemize}
\end{defn}
\par
The matrices $G$, $G'$ and $G''$ are defined in the following way. $G_{Ii}$ (similarly $G'_{Ii}$,$G''_{Ii}$) is the number of times an edge of a tetrahedron $\Delta_i$ with parameter $Z_i$, (similarly $Z'_i$, $Z''_i$) is glued to edge $E_I$. Thus for each edge $E_I$, we have the following gluing equation: 
$$E_I:\quad \sum_{i=1}^{N}(G_{Ii}Z_i+G'_{Ii}Z'_i+G''_{Ii}Z''_i)=2\pi i \quad I=1,...,N.$$
This additive form of the gluing equations is the log of the multiplicative gluing equations, with $Z_i:=log(z_i)$ (similarly for $Z'_i$, $Z''_i$).
\par
\begin{prop}\label{propsqua}If $\mathcal{T}$ has $N$ tetrahedra, then it has $N$ edges as well.
\end{prop}
\begin{proof}
	
	This follows from an Euler characteristic argument for the boundary torus. Let $f, e, v$ be the number of faces, edges, and vertices of the triangulation of the boundary torus, respectively, and let $N_t, N_e$ be the number of tetrahedra, edges of $\mathcal{T}$. Each edge of $\mathcal{T}$ meets the boundary torus twice, and each intersection corresponds to a vertex on the boundary torus, so $v=2N_e$. Each tetrahedron of $\mathcal{T}$ meets the boundary torus 4 times, and each intersection corresponds to a face on the boundary torus, so $f=4N_t$. Finally, each edge of the triangulation of the boundary torus corresponds to two glued faces of tetrahedra, and each face of each tetrahedron meets the boundary torus 3 times, so $e=\frac{3*4*N_t}{2}=6N_t$. Since the Euler characteristic of the torus is 0, $f-e+v=4N_t-6N_t+2N_e=0$, so $N_t=N_e$.
\end{proof}
Thus  $G$, $G'$ and $G''$ are square matrices.
\par
A \textit{flattening} $\bm{f}$ is a triple $f,f',f''\in \mathbb{Z}^N$ such that 
$$f+f'+f''=[1,...,1]^t,$$ $$Gf+G'f'+G''f''=[2,...,2]^t.$$
\par
The matrices $G_{\mu}$, $G'_{\mu}$ and $G''_{\mu}$ are defined by deleting a row of $G$, $G'$ and $G''$, and replacing them with the information for the completeness equation of a meridian $\mu$ of the knot:
$$\mu:\quad \sum_{i=1}^{N}(C_iZ_i+C'_iZ'_i+C''_iZ''_i)=0.$$
The gluing equations have redundancy (see \cite{dimofte2013quantum}), so we can remove a row without losing any information.
\par
If our triangulation $\mathcal{T}$ is essential, then there exists a geometric solution to the gluing equations and the completeness equation for the meridian. As such we will fix a solution $z^0=(z_1^0,...,z_N^0)$ to define the 1-loop invariant:
\begin{defn}[Dimofte-Garoufalidis \cite{dimofte2013quantum}]
	$$\tau(\mathcal{T}):=\pm \frac{1}{2}\frac{det[G_\mu diag(\zeta_1,...,\zeta_N)+G'_\mu diag(\zeta'_1,...,\zeta'_N)+G''_\mu diag(\zeta''_1,...,\zeta''_N)]}{\prod_{i=1}^{N}\zeta_i^{f_j}\zeta_i^{'f'_j}\zeta_i^{''f''_j}}$$
	where $$\zeta_j:=\frac{1}{z_j}, \quad \zeta'_j:=\frac{1}{(1-z_j)}, \quad \zeta''_j:=\frac{1}{z_j(z_j-1)},$$ and this is evaluated at $z^0$.
\end{defn}
It was proved in \cite{dimofte2013quantum} that this invariant is independent of choice of \textit{quad type} for each tetrahedron (a choice of which pair of edges to label $z$ as opposed to $z'$ or $z''$), choice of edge equation to remove, choice of path representing the meridian of the knot, and choice of flattening. Additionally, it was shown to be invariant under 2-3 Pachner moves, and later additionally under 0-2 Pachner moves \cite{pandey2023geometry}, and since the set of essential ideal triangulations is connected under 2-3 and 0-2 Pachner moves \cite{kalelkar2024connecting}, this is independent of triangulation and a manifold invariant.

\subsection{The Ohtsuki-Takata Invariant}
\begin{figure}[H]
	\centering
	\includegraphics[width=.75\linewidth]{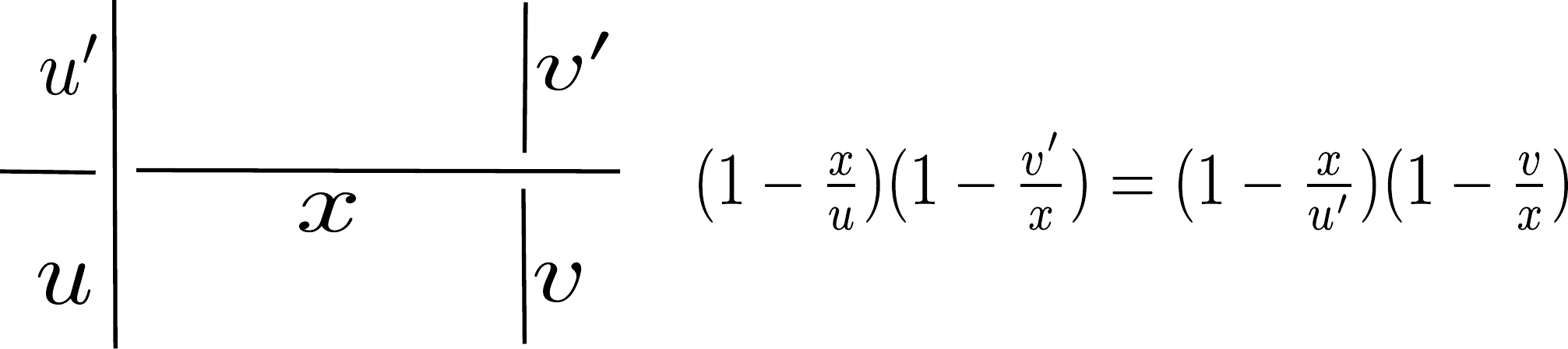}
	\caption{Alternating Hyperbolicity Condition}
	\label{fig:HypEq}
\end{figure}

\begin{figure}[H]
	\centering
	\includegraphics[width=.6\linewidth]{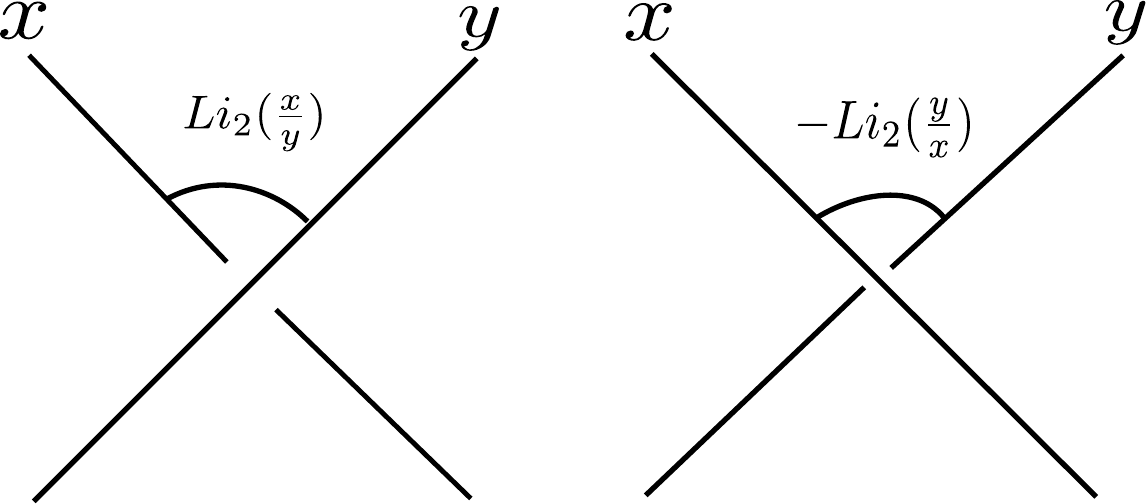}
	\caption{Corner Contributions to $V$}
	\label{fig:DiLog}
\end{figure}

\begin{figure}[H]
	\centering
	\includegraphics[width=.65\linewidth]{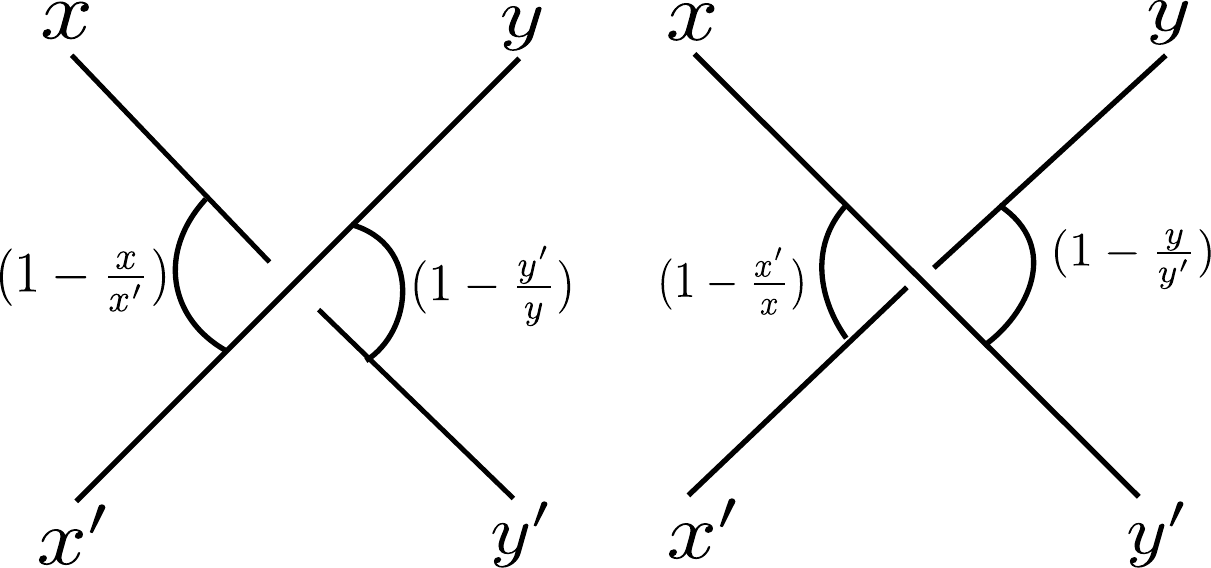}
	\caption{Crossing Contributions to $\Omega_1$}
	\label{fig:Omega1Cross}
\end{figure}

\begin{figure}[H]
\centering
\includegraphics[width=.65\linewidth]{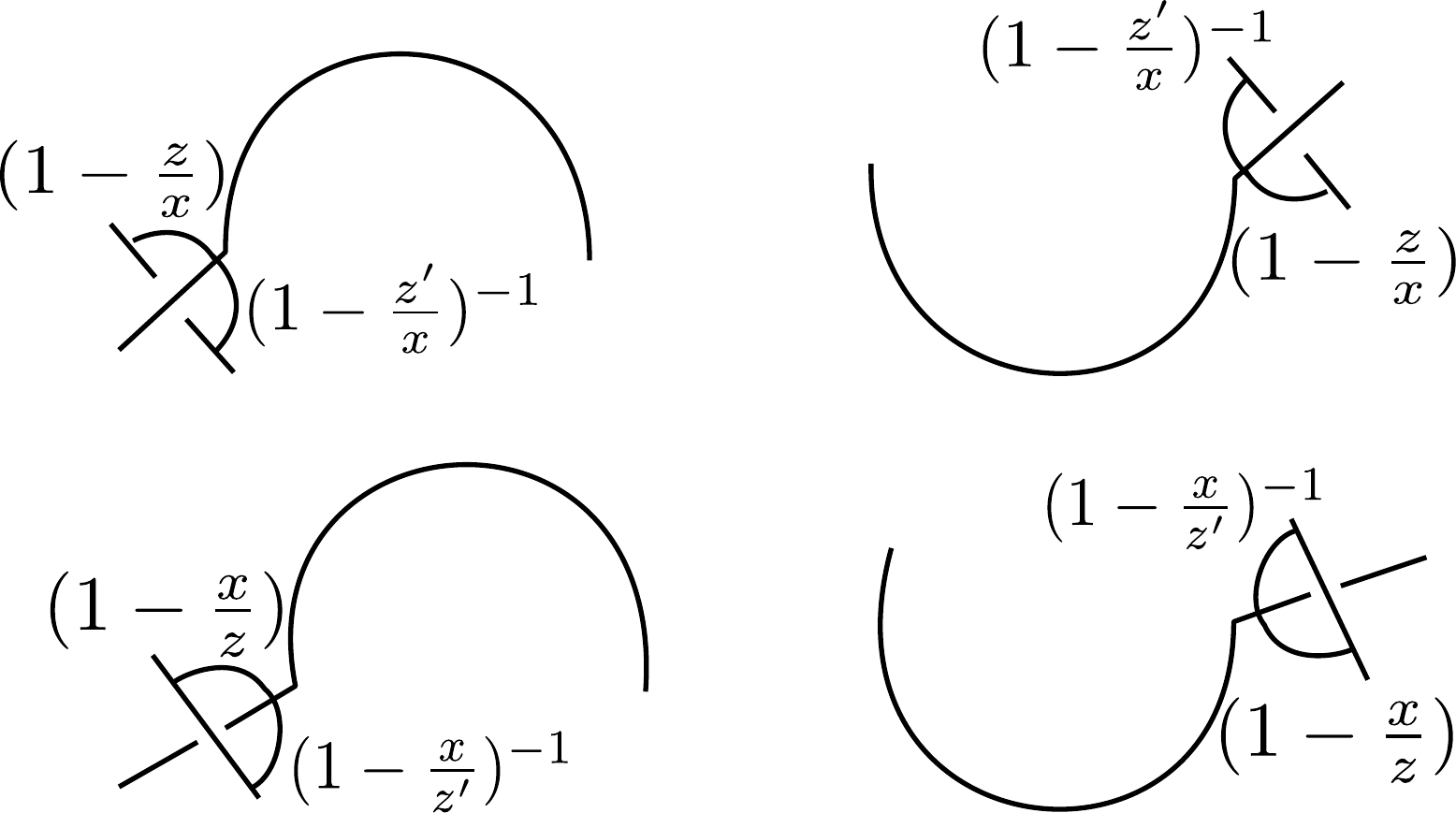}
\caption{Cup/Cap Contributions to $\Omega_1$}
\label{fig:Omega1Caup}
\end{figure}

\begin{figure}[H]
	\centering
	\includegraphics[width=.65\linewidth]{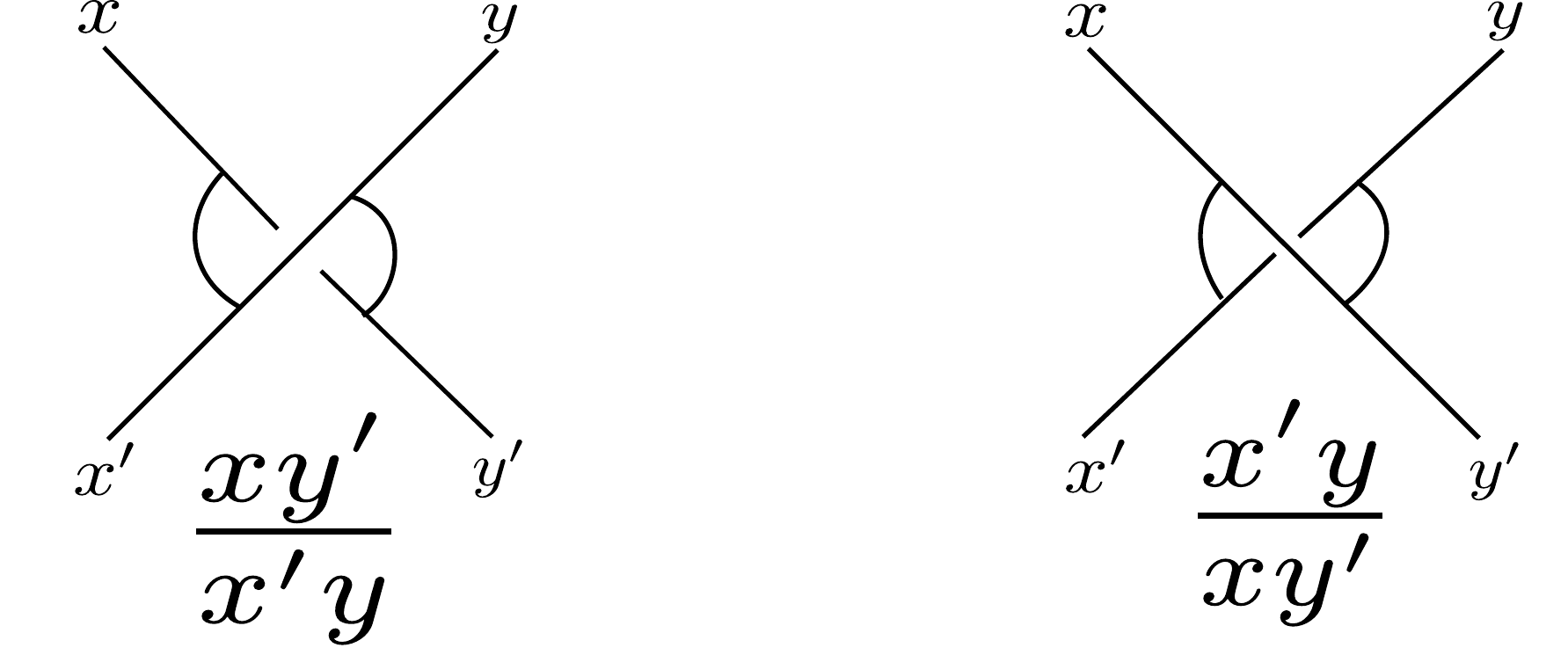}
	\caption{Crossing Contributions to $\Omega_2$}
	\label{fig:Omega2}
\end{figure}

\begin{figure}[H]
	\centering
	\includegraphics[width=.5\textwidth]{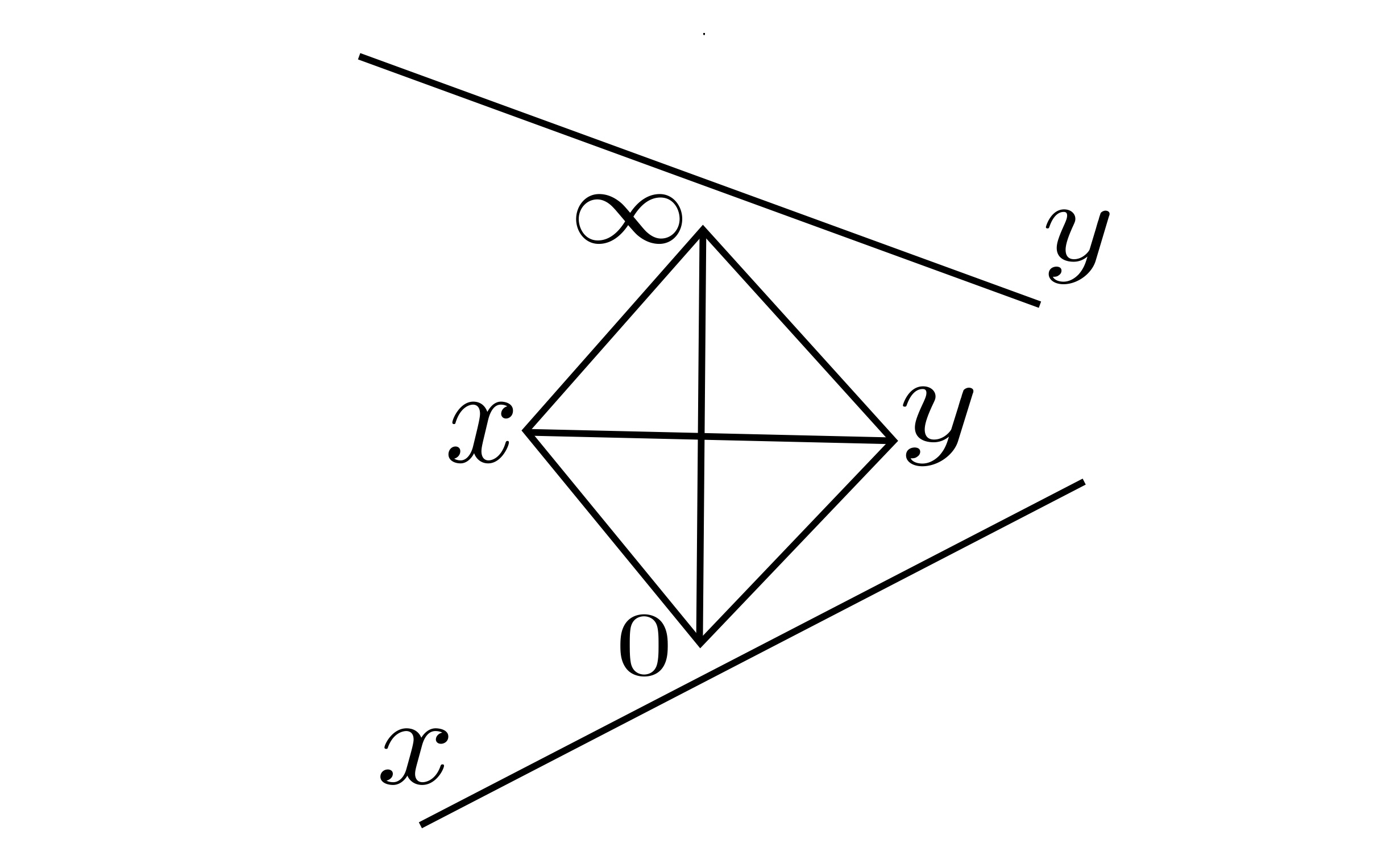}
	\caption{A Tetrahedron with Shape Parameter $\frac{y}{x}$}
	\label{fig:crosstet}
\end{figure}

The Ohtsuki-Takata invariant is defined for an open knot diagram, where the open ends of the diagram are along an alternating segment. The segments of the diagram receive labels as in \cite{yokota2003jones}. The open segments don't receive a label, and any segments connected to an open segment by a series of overpasses are labeled $0$, while any segments connected to an open segment by a series of underpasses are labeled $\infty$. In the case of alternating knots, this means there is one segment labeled $0$, and one labeled $\infty$. All segments adjacent to the unbounded region are labeled 1, and all remaining segments are labeled with non-zero complex parameters $x_1,...,x_n$. These parameters are required to satisfy the \textit{hyperbolicity equations}, which in the alternating case is given as in Figure \ref{fig:HypEq}. 
For an open diagram with such labels, we say that a corner is an \textit{essential corner} if the ratio of the labels of the segments meeting at that corner is finite and non-constant. That is, if the corner is not in the unbounded region, and does not meet an $\infty$-segment or a $0$-segment. To define the \textit{Yokota potential function} $V(\mathring{D})$, a signed dilogarithm is associated to each essential corner as shown in Figure \ref{fig:DiLog}, and the value of the potential function is the sum of these over all essential corners.
Additionally, we define two normalizing functions associated with $\mathring{D}$, $\Omega_1$ and $\Omega_2$. $\Omega_1$ is defined for corners at crossings and cups/caps as shown in Figures \ref{fig:Omega1Cross} and \ref{fig:Omega1Caup}, and defined for $\mathring{D}$ to be the product over all of the values at crossings, cups, and caps where the corners chosen are essential.

In \cite{ohtsuki2015kashaev}, Ohtsuki and Takata defined $\Omega_2$ in terms of an orientation of $\mathring{D}$, but Garoufalidis and Yoon (\cite{garoufalidis20261}, proof of lemma 4.9) showed that this is equivalent to a formulation that does not involve an orientation, so this is what we will describe here. $\Omega_2$ is defined for corners at each crossing as in Figure \ref{fig:Omega2}, and defined for $\mathring{D}$ as the product of this value at essential corners over all crossings. Note that if in Figure \ref{fig:Omega2}, the left (resp. right) is not essential, we do not include the contribution of $(x'/x)^{\pm 1}$ (resp. $(y'/y)^{\pm 1}$).
Define the critical point equations of $V(\mathring{D})$ to be 
$$exp(x_i\frac{\partial V}{\partial x_i})=1 \quad i=1,...,n.$$ A \textit{critical point} of $V$ is a solution to these equations. As shown in \cite[Theorem 4.1]{murakami2018volume}, there is a correspondence between solutions to the critical point equations, and solutions to the gluing and completeness equations. Under this correspondence, essential corners of $\mathring{D}$ correspond to tetrahedra of $\mathcal{T}(\mathring{D})$ and the shape parameters of a tetrahedron is given by the ratio of the Yokota labels of the segments at that corner. See Figure \ref{fig:crosstet}. 
As such there exists a critical point of $V(\mathring{D})$ which corresponds to the geometric solution to the gluing and completeness equations. We will call this point $x^0=(x^0_1,...,x^0_n)$.
\begin{defn}[Ohstuki-Takata\cite{ohtsuki2015kashaev}] The \textit{Ohtsuki-Takata invariant} of an essential open diagram $\mathring{D}$ is $$\omega(\mathring{D}):=\pm \frac{\Omega_1 \Omega_2}{2}det(x_i\frac{\partial}{\partial x_i}(x_j\frac{\partial V}{\partial x_j})), \quad i,j\in \{1,...,n\},$$ where this is evaluated at $x^0$.
\end{defn}

\section{Alternating Braids}
Both the Ohtsuki-Takata and 1-loop invariants are expressed in terms of the determinant of a matrix corresponding to some system of equations. The main result of this paper is to relate these two matrices in the case of reduced alternating braid diagrams.
\begin{thm}\label{mainthm}
	Let $K$ be a hyperbolic knot, let $\mathring{D}$ be an open, reduced (as a knot diagram), alternating braid diagram of $K$, and let $\mathcal{T}$ be the triangulation associated to the collapsed octahedral decomposition coming from $\mathring{D}$. Then for the values defined as in $\omega(\mathring{D})$ and $\tau(\mathcal{T})$, $$det(x_i\frac{\partial}{\partial x_i}(x_j\frac{\partial V}{\partial x_j}))=\pm (\prod_{k=1}^{N}z_k) det(G_{\mu}diag(\zeta)+G'_{\mu}diag(\zeta')+G''_{\mu}diag(\zeta'')).$$
\end{thm}
\textbf{Note:} To show the full equivalence between the 1-loop and Ohstuki-Takata invariants for this case, it remains to prove the following normalization identity: $$\Omega_1\Omega_2(\prod_{k=1}^{N}\zeta_j^{f_j}(\zeta'_j)^{f'_j}(\zeta''_j)^{f''_j})(\prod_{k=1}^{N}z_k)=\pm 1$$

\subsection{Labeling of Yokota Variables}
\textbf{Note:} In this section, by a diagram of ``\textbf{chosen type}" we will mean a reduced (as a knot diagram), alternating braid diagram of a hyperbolic knot.
\par
We will be considering the open knot diagram that is the closure of all of the strands except for the first one. Without loss of generality, consider the case where leftmost crossings are negative crossings, $\sigma_1^{-1}$. Since we are considering only alternating braids, we will have only positive crossings between even numbered strands on the left and odd numbered strands on the right, and and only negative crossings between odd numbered strands on the left and even numbered strands on the right. Essentially, our braid word will only consist of $\sigma_i$ for $i$ even and $\sigma_i^{-1}$ for $i$ odd. 
\par
\begin{lem}\label{lem1}For any diagram $D$ of \textbf{chosen type}, and any $2\leq i\leq m-1$, where $m$ denotes braid width, $D$ is equivalent up to conjugation to a braid whose word is expressed as $...\sigma_i^{\pm 1}...\sigma_{i-1}^{\mp 1}...\sigma_i^{\pm 1}...\sigma_{i-1}^{\mp 1}$.
\end{lem}
\begin{proof}
Assume the braid word of D has only one $\sigma_{i-1}^{\mp 1}$. (In this proof I will also refer to the crossings represented by $\sigma_i$ as $\sigma_i$ for the sake of brevity.) Then we can draw a circle meeting that crossing and continuing vertically and passing through the arcs that close our braid, so as to only meet the diagram at that crossing. This shows that this is a reducible crossing. Thus D is not reduced, contradiction. Therefore we must have at least two $\sigma_{i-1}^{\mp 1}$'s. 
\par
Assume now that all copies of $\sigma_{i-1}^{\mp 1}$ are consecutive. Then cutting the diagram vertically between the $\sigma_{i-1}^{\mp 1}$'s and the $\sigma_i^{\pm 1}$'s decomposes the braid into a  connect sum of two braids. Both components have only positive crossings in even columns and negative crossings in odd columns, so their reduced words are nonempty. Therefore both components are nontrivial, and the braid is a compound knot, and thus not hyperbolic, and we have another contradiction. Thus we must have that not all copies of $\sigma_{i-1}^{\mp 1}$ are consecutive. 
\par
Furthermore, the braid cannot be equivalent to any braid in which all of the $\sigma_{i-1}^{\mp 1}$'s are consecutive. As such we cannot have that all of the $\sigma_{i-1}^{\mp 1}$'s are at the beginning and end of the braid word, as conjugation by some number of copies of $\sigma_{i-1}^{\mp 1}$ would give that this braid is equivalent to one who has all copies of $\sigma_{i-1}^{\mp 1}$ consecutive. Additionally, it cannot be that the only crossings between the $\sigma_{i-1}^{\mp 1}$'s are copies of $\sigma_j^{\pm1}$ for $|i-1-j|\geq 2$, as we could use the relation that $\sigma_i^{\pm1}\sigma_j^{\pm1} \sim \sigma_j^{\pm1}\sigma_i^{\pm1}$ for $|i-j|\geq 2$ to show once again that this braid is equivalent to one whose copies of $\sigma_{i-1}^{\mp 1}$ are all consecutive. Thus up to conjugation we must have $...\sigma_i^{\pm 1}...\sigma_{i-1}^{\mp 1}...\sigma_i^{\pm 1}...\sigma_{i-1}^{\mp 1}$.     
\end{proof} 
\begin{rmk}
	Conjugation of a braid preserves its closure, and therefore does not change the knot. Thus we will consider our braid diagrams to have the form $...\sigma_1^{-1}...\sigma_2...\sigma_1^{- 1}...\sigma_2$.
\end{rmk}
For a knot diagram $D$ with Yokota labels, let $\mathcal{G}(D)$ be the graph whose vertices correspond to segments of D with assigned Yokota parameters, and whose edges correspond to essential corners, where the vertices at the end of each edge correspond to the segments separated by the essential corner. Let $\mathcal{G}'(D)$ be the subgraph whose vertices are only those with non-constant Yokota parameters. 
\par
\begin{lem}\label{lem2}For a diagram $D$ of \textbf{chosen type}, there exists a Hamiltonian path in $\mathcal{G}'(D)$ which begins at a vertex adjacent in $\mathcal{G}(D)$ to a vertex corresponding to a segment adjacent to the unbounded region (a segment labeled 1).
\end{lem}
\par
The proof of this lemma is done by construction for four cases, corresponding to four types of behavior around the $\infty$-segment. The $\infty$-segment is always adjacent at the top(start) to the first $\sigma_1^{-1}$ crossing of the braid word, but it can be adjacent at the bottom(end) to either a $\sigma_1^{-1}$ crossing or a $\sigma_2$ crossing. Therefore the behavior around this segment can be described in terms of the portion of the braid word around the first $\sigma_1^{-1}$ crossing:
\begin{enumerate}
	\item The first two $\sigma_1^{-1}$'s in the braid word are consecutive up to braid relations, which is to say there is no $\sigma_2$ between them.
	\item The first two $\sigma_1^{-1}$'s in the braid word are \textbf{not} consecutive up to braid relations. There is no $\sigma_3^{-1}$ between the $\sigma_2$ preceding the first $\sigma_1^{-1}$ and the $\sigma_2$ following it.
	\item The first two $\sigma_1^{-1}$'s in the braid word are \textbf{not} consecutive up to braid relations. There is a positive odd number of $\sigma_3^{-1}$ between the $\sigma_2$ preceding the first $\sigma_1^{-1}$ and the $\sigma_2$ following it.
	\item The first two $\sigma_1^{-1}$'s in the braid word are \textbf{not} consecutive up to braid relations. There is a positive even number of $\sigma_3^{-1}$ between the $\sigma_2$ preceding the first $\sigma_1^{-1}$ and the $\sigma_2$ following it.
\end{enumerate}
Examples of the first 3 cases are shown in Figure \ref{fig:startbehavior}, with the fourth case being akin to the third, but with an even number of $\sigma_3^{-1}$ crossings in place of an odd number.

\begin{proof}Each segment in our diagram exists between two columns of crossings, so we will say it ``exists" in a certain column. For example: segments to the left of any crossings will be in column 1, segments between $\sigma_1^{-1}$ and $\sigma_2$ crossings will be in column 2, and so on. We will also refer to this Hamiltonian path as traversing segments, when really what is meant is traversing the vertices corresponding to the segments, but as these are identified, this poses no problems.
\par
\begin{figure}[H]
    \centering   
\begin{center}
\begin{tikzpicture}
\pic[
rotate=0,
braid/.cd,
every strand/.style={ultra thick},
strand 1/.style={black},
strand 2/.style={black},
strand 3/.style={black},
strand 4/.style={black},
] {braid={s_2 s_1^{-1} s_1^{-1}}};
\node[label={45:$\infty$}] at (.8,-2.4) {};
\node[label={45:$1$}] at (-.5,-2.6) {};
\node[label={45:$1$}] at (-.5,-3.7) {};
\node[label={45:$1$}] at (.5,-.6) {};
\node[label={45:$1$}] at (.5,-1.6) {};
\end{tikzpicture}
\quad\quad\quad\quad
\begin{tikzpicture}
\pic[
rotate=0,
braid/.cd,
every strand/.style={ultra thick},
strand 1/.style={black},
strand 2/.style={black},
strand 3/.style={black},
] {braid={s_2 s_1^{-1} s_2}};
\node[label={45:$\infty$}] at (.8,-2.4) {};
\node[label={45:$1$}] at (-.5,-2.7) {};
\node[label={45:$1$}] at (.5,-.6) {};
\node[label={45:$1$}] at (.5,-1.6) {};
\end{tikzpicture}
\quad\quad\quad\quad
\begin{tikzpicture}
\pic[
rotate=0,
braid/.cd,
every strand/.style={ultra thick},
strand 1/.style={black},
strand 2/.style={black},
strand 3/.style={black},
strand 4/.style={black},
] {braid={s_2 s_1^{-1} s_3^{-1} s_2}};
\node[label={45:$\infty$}] at (.8,-2.8) {};
\node[label={45:$1$}] at (-.5,-2.9) {};
\node[label={45:$1$}] at (.5,-.6) {};
\node[label={45:$1$}] at (.5,-1.6) {};
\end{tikzpicture}
\end{center}
\caption{First, second and third cases for behavior around $\infty$-segment}
\label{fig:startbehavior}
\end{figure}
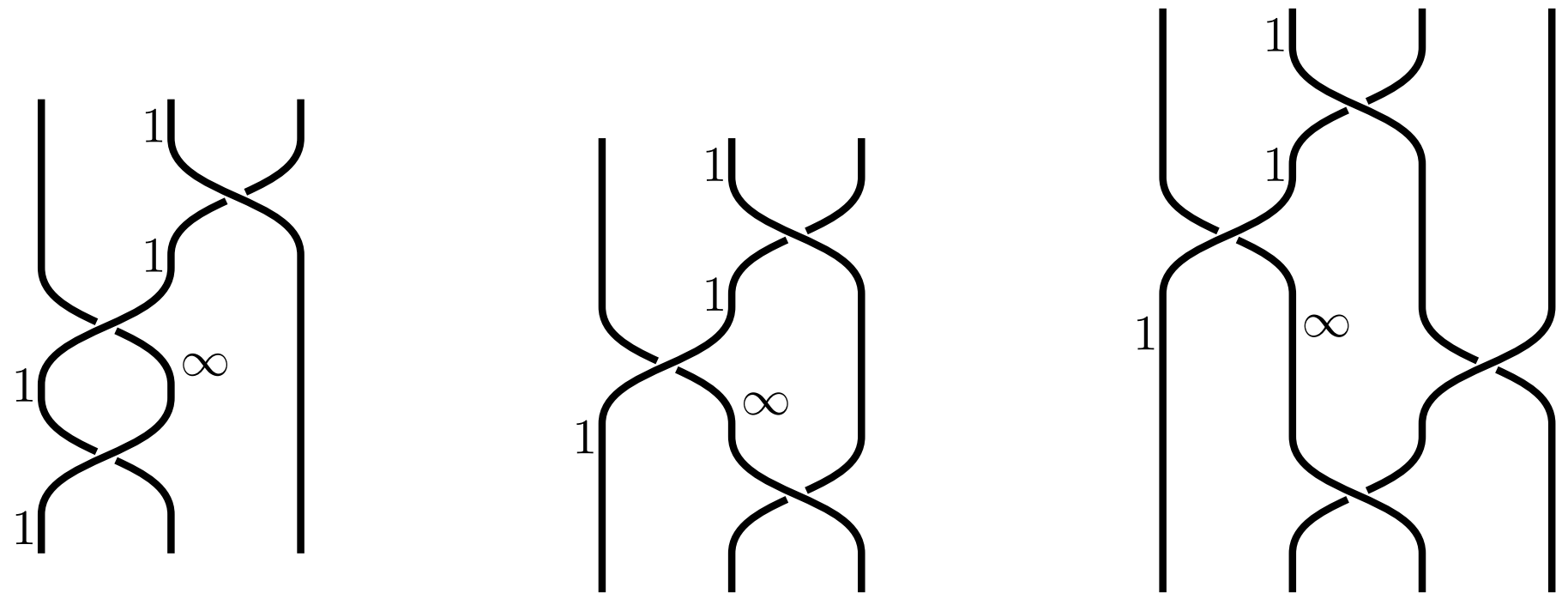
Note that every segment is adjacent to the segment above and below it (assuming these segments have non-constant labels), because these segments share a corner, which must be essential because both segments have non-constant labels. In column 1 none of the segments have non-constant labels, in column 2 the $\infty$-segment and $0$-segment separate all of the non-constant segments from the segments with label 1, and in the column 3 or greater, all of the segments have non-constant labels, with the top and bottom of the braid being connected by a braid closure, so we can consider these segments to be ``above" the segment just below the cap and also ``below" the segment just above the cup.
\par 
\textbf{Case 1:} The path starts on the first segment in column 2 below the $\infty$-segment. This segment is separated from a segment in column 1 adjacent to the unbounded region by an essential corner at its top crossing's south corner. The path then proceeds along successive segments down column 2 until one of those segments terminates at a $\sigma_2$ crossing. Then this segment is adjacent to a segment in column 3, and the path proceeds to that segment. The path then continues by proceeding along successive segments \textbf{up} column 3 until one of those segments connects to a $\sigma_3^{-1}$ crossing. This pattern then continues, alternating going down even columns and up odd columns, until the path reaches the rightmost column. Then the path continues all the way around the rightmost column (either up or down depending on the parity of braid width), then when it has passed through every segment in that column it passes the next column to the left. This passage is along the lower(resp. upper) two strands at the crossing at which the path passed to the rightmost column along the upper(rep. lower) two strands. The path then passes along the second rightmost column in the same direction it was traversed on the way in, passing through all of the remaining segments, then passing back left to the next column. This pattern continues along the columns from right to left until the path terminates at the final non-constant segment, which is in column 2, except in the case where the $0$-segment meets the first $\sigma_2$ crossing after the initial $\sigma_1^{-1}$ crossings, in which case it is in column 3.
\par
\textbf{Case 2:} The path starts at the segment in column 3 to the right of the $\infty$-segment, adjacent at their lower crossing. This is separated from a segment labeled 1 by the south corner of its upper crossing. The path proceeds up column 3 then continues in the same way as case 1, until its termination.
\par
\textbf{Case 3:} The path starts at the segment in column 3 that is directly above the final $\sigma_2$ crossing before the first $\sigma_1^{-1}$ crossing. This is separated from a segment labeled 1 by the north corner of this $\sigma_2$ crossing. This case proceeds the same way as the previous cases, until the path is passing along columns to the left and going down along column 4. When it reaches the segment adjacent to the start segment, the path passes down one more segment, then left to the segment in column 3, which is adjacent by the existence of a $\sigma_3^{-1}$ crossing. The path then passes down a segment, then to the right along the lower strands of the same $\sigma_3^{-1}$ crossing. This pattern continues through to the lower right strand of the last of the consecutive $\sigma_3^{-1}$ crossings, and then continues down to the last segment of column 4, passes to column 3, and terminates as before. 
\par
\textbf{Case 4:} The path starts at the segment in column 3 that is directly below the final $\sigma_2$ crossing before the first $\sigma_1^{-1}$ crossing. This is separated from a segment labeled 1 by the south corner of this $\sigma_2$ crossing. This case then proceeds in the same way as case 3.
\par
In each of these cases, the path passes through every non-constant segment in every column.
\end{proof}

\par
\begin{figure}[H]
        \centering
        \includegraphics[width=.3\textwidth]{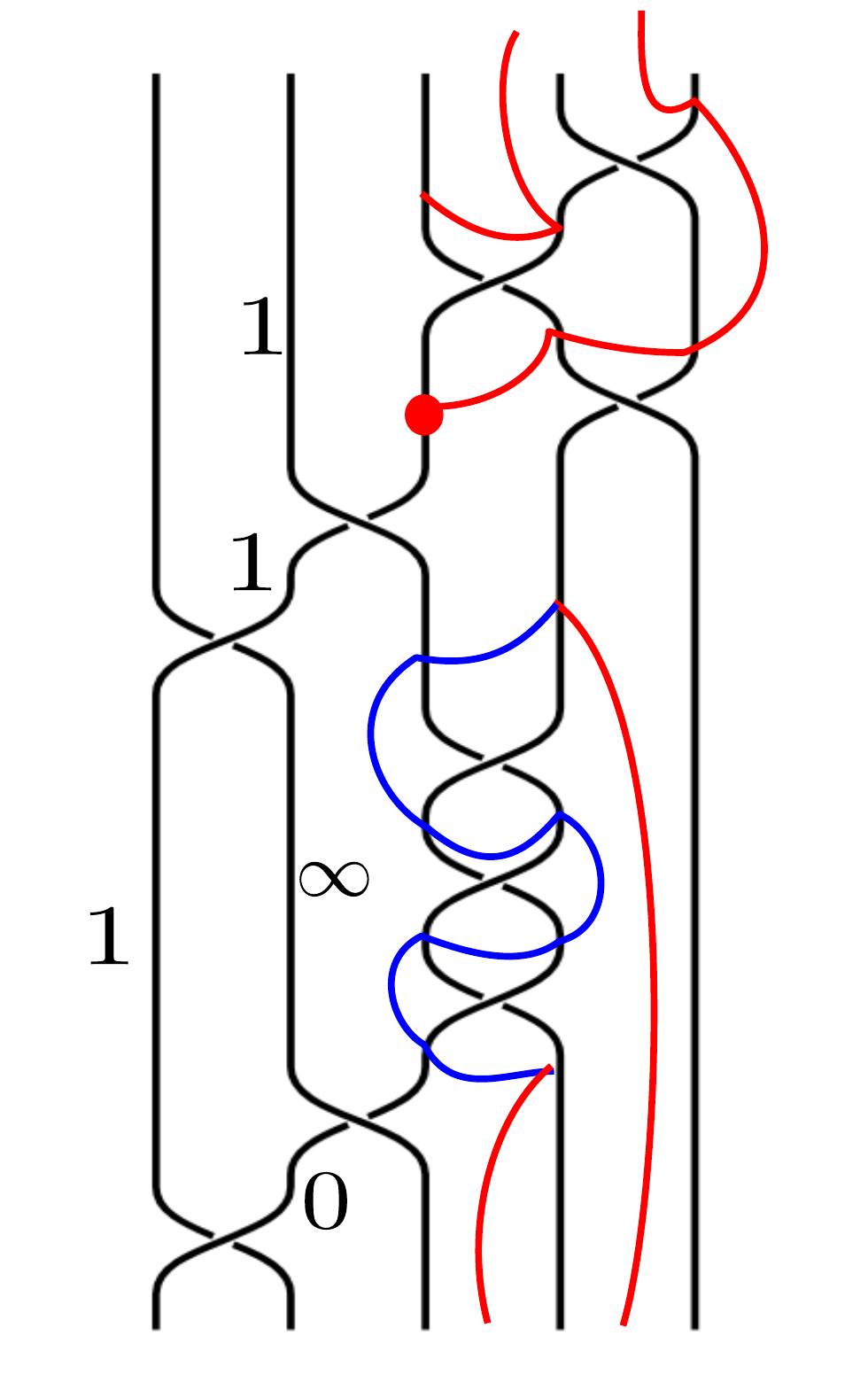}
        \caption{A Hamiltonian Path on a Braid Diagram}
        \label{fig:hampath}
    \end{figure}

 To give an example, see Figure \ref{fig:hampath}, where the Hamiltonian path is redrawn on top of the braid diagram, such that a line segment between segments of the braid corresponds to an edge of the Hamiltonian path between associated vertices. This braid is shown as entirely unclosed, but we consider the closure of every strand except the leftmost. This is an example of the third case. The Hamiltonian path is shown in red, with the dot representing the start, and the deviation from the usual pattern needed in the third case shown in blue.

\par This path will give us our labels for the parameters of the Yokota potential function. The first segment will receive the label $x_1$, and subsequent segments along the path will receive the labels $x_2,...,x_n$. In the next section we will discuss a corresponding elimination of shape parameters.

\subsection{Shape Parameter Elimination}

\begin{lem}[\cite{garoufalidis20261}, Lemma 4.7]\label{lem_elim}
	Let $f_1,...,f_k$ be functions in variables $x_1,..,x_k$ and fix a solution to $$exp(f_1)=...=exp(f_k)=1.$$ Suppose that $exp(f_k)$ is of the form $g(x_1,...,x_{k-1})x_k$ so that we can substitute $x_k$ with $1/g(x_1,...,x_{k-1})$ by solving $exp(f_k)=1$. then we have $$det(\frac{\partial (f_1,...,f_k)}{\partial (x_1,...,x_k)})=\frac{1}{x_k}det(\frac{\partial (f'_1,...,f'_{k-1})}{\partial (x_1,...,x_{k-1})})$$ where $f'_i$ is obtained from $f_i$ by substituting $x_k$ with $1/g(x_1,...,x_{k-1})$.
\end{lem}

We know that gluing equations for C and R-edges have the form $z_{i_1}...z_{i_n}=1$, so we can use them to apply Lemma \ref{lem_elim} and reduce the size of the system. Specifically, they can be used to eliminate any shape parameter corresponding to a tetrahedron at the crossing for a C-edge, or in that region for an R-edge. Now we will show how to apply this to eliminate every shape parameter, except for one for each Yokota parameter $x_i$, which corresponds to the tetrahedron which lies at the essential corner between segments with labels $x_i$ and $x_{i-1}$, or $1$ and $x_1$. To do so we will construct a map that gives a choice of R-edge or C-edge for each corner to be eliminated, save for one, which will be eliminated by the completeness equation for the meridian.
\par
\textbf{Note:} As before, we will refer to segments as lying in columns, or ``segment columns", but now we will also refer to the columns in which crossings lie as ``crossing columns". Specifically, crossing column $i$ will lie between segment columns $i$ and $i+1$. As all bounded regions (except for the central region) are bounded above and below by crossings in the same column, we will refer to those regions as lying in crossing columns as well. Additionally, not including behavior around the $\infty$-segment for braids of the third and fourth type, there is one crossing in the column $i$ for $2 \leq i\leq m-1$ in which our previously chosen path passes from segment column $i$ to $i+1$. In the 2nd column, this may be where the path begins, terminates, or both. We will refer to this as the \textbf{special crossing} of that column. 
\par
\begin{lem}\label{lem3}
    For every diagram of \textbf{chosen type}, and labeling determined as previously, there exists a map $\kappa$ from the set of crossings with corresponding C-edges and regions with corresponding R-edges to the set of all essential corners such that the complement of the image of the map is one corner for each $x_i$, lying between segments with labels $x_i$ and $x_{i-1}$, or $1$ and $x_1$, (possibly) union one other corner.
\end{lem}
The proof is a construction, with general behavior for columns 4 and on, and the 4 cases for behavior around the $\infty$-segment. We will refer to corners of crossings with north, south, east, and west with respect to the previously used orientation of the braid diagram.
\begin{proof}
    \textbf{General behavior:} The central (rightmost) region maps to the east corner of the rightmost special crossing. Consider the $i$th column, for $4\leq i \leq m-1$. The $i$th column special crossing maps to its west corner. All regions above this crossing and below the region adjacent to the ($i$-1)th special crossing map to the south corner of the crossing above it, and all crossings above the $i$th special crossing and below this region map to their north corner. All regions below the $i$th special crossing and above the region adjacent to the ($i$-1)th special crossing map to the north corner of the crossing above it, and all crossings below the $i$th special crossing and above this region map to their south corner. The region adjacent to the ($i$-1)th special crossing maps to the east corner of the ($i$-1)th special crossing. Therefore in the $i$th column, the $\kappa$ maps onto the north and south corner of every crossing except the special crossing, where it maps onto the east and west corners. This means that what remains is the east and west corner for everything except the special crossing, and the north and south corner of the special crossing. This is exactly the corners that lie between consecutive segments on our chosen path.
    \par
 	\textbf{Case 1:} The crossings and regions in column 3 map in the same way as in the general behavior. Two regions in column 2, adjacent to the $\infty$ and $0$-segments, do not have R-edges. Above the region adjacent to the $\infty$-segment and below the region adjacent to the $0$-segment, crossings map to their north corners and regions map to the south corner of the crossing above them. In this portion of the braid there are no essential corners at west corners, so what remains is the east corners, and the south corner of the crossing above the region adjacent to the $\infty$-segment. Above the region adjacent to the $0$-segment and below the special crossing, crossings map to their south corner and regions map to the north corner of the crossing below them. The special crossing has a crossing edge if and only if its west corner is essential, in which case it maps to that corner. In this portion of the braid, what remains is the east and west corners of non-special crossings, as well as the north and potentially south corner of the special crossing. In column 1, all crossings and regions with C and R-edges lie between two regions with no R-edges. Crossings map to their south corners, and regions map to the north corner of the crossing below them. There are no west essential corners in this column. Therefore the remaining corners not in the image of $\kappa$ in these columns are exactly the ones which lie between consecutive segments, with the exception of the corner separating the start of the path from a segment with label 1, and the south corner of the crossing in column 2 above the region adjacent to the $\infty$-segment.
    \par
    \textbf{Case 2:} The crossings and regions in column 3 map in the same way as in the general behavior. The crossings and regions in columns 1 and 2 map in the same way as in case 1, with the caveat that the special crossing in column 2 does not have a crossing edge, or an essential corner at its north corner. The remaining corners not in the image of $\kappa$ in these columns are exactly the ones which lie between consecutive segments, with the exception of corner separating the start of the path from a segment with label 1, and possibly the north corner of the crossing in column 1 below the region adjacent to the $\infty$-segment, in the case that there are more than two crossings in column 1.
    \par
    \textbf{Case 3:} The special crossing in column 2 must be below the crossings in column 3 around which snakes the path and above the special crossing in column 3, because the path traverses segment column 3 downwards. In column 3, the special crossing maps to its west corner, and above the special crossing and below the region adjacent the special crossing in column 2, regions map to the south corner of the crossing above it, and crossings map to their north corner. The region adjacent to the special crossing in column 2 maps to its east corner. The crossings above this region are the ones around which the path snakes. Regions in between these crossings map to the north corner of the crossing below it, and the crossing maps alternate mapping to the east and west corners, whichever side does not separate adjacent segments along the path. The region above these crossings maps to the east corner of the lowest crossing in column 2 to which it is adjacent. Above this region and below the special crossing, crossings map to their south corner and regions map to the north corner of the crossing below them. In column 2, this is the same as in the previous cases, except that the crossing above the region adjacent to the $\infty$-segment maps to its south corner. Column 1 is the same as in the previous cases. The remaining corners not in the image of $\kappa$ in these columns are exactly the ones which lie between consecutive segments, with the exception of corner separating the start of the path from a segment with label 1, and possibly the north corner of the crossing in column 1 below the region adjacent to the $\infty$-segment, in the case that there are more than two crossings in column 1.
    \par
    \textbf{Case 4:} This case is exactly the same as Case 3, with the exception that the region adjacent to the start segment in column 3 maps to the north corner of the crossing below it, and the crossing in column 2 above the region adjacent to the $\infty$-segment maps to its north corner rather than its south corner.
    \par
    This map is well-defined, as we did not double-count any crossing or region, and this coincides with the argument in \cite{garoufalidis20261} about the number of variables eliminated.
\end{proof}
\subsection{Completeness Equations for Alternating Braids}
In this section, we will discuss how to find a completeness equation of the meridian for any alternating braid diagram, and how to use it to eliminate the remaining shape parameter.
\par
How to find the meridian completeness equation is slightly different between the first case and the second, third, and fourth cases. The first case is similar to the examples in \cite{garoufalidis20261} and \cite{murakami2018volume}, so we will explain in detail what this looks like in the latter cases.
\par
\begin{prop}\label{prop1}
	Assume that for a diagram $\mathring{D}$ of \textbf{chosen type}, there exists an essential corner that is neither in the image of $\kappa$, or the unique corner between segments with labels $x_{i-1}$ and $x_i$ for some $i$. Let this corner have corresponding shape parameter $z_j$. Then the completeness equation of a meridian of the knot can be written in the form $g(z_1,...,z_{j-1},z_{j+1},...,z_N)z_j=1$.
\end{prop}
\begin{proof}
Let $\mathcal{V}(+\infty)$, $\mathcal{V}(-\infty)$, and $\mathcal{V}(K)$ represent the un-collapsed triangulations of boundaries of neighborhoods of the point at $+\infty$, the point at $-\infty$, and the knot, respectively, and let $\mathcal{V'}(+\infty)$, $\mathcal{V'}(-\infty)$, and $\mathcal{V'}(K)$ represent the corresponding collapsed triangulations. In this section it may be helpful to refer to \cite{kim2018octahedral}, which shows twisted octahedra, and helps show how they glue. The faces at truncated vertices form the aforementioned triangulations. 
\par
Each octahedron contributes to $\mathcal{V}(+\infty)$ a "pointed bow tie", where the inner edges of the bow tie correspond to upper faces of the octahedron, and outer edges of the bow tie correspond to lower faces of the octahedron, see Figure \ref{fig:bowtie}. Additionally, vertices of the bow tie correspond to edges of the octahedron, with the vertex at the center corresponding to the identified upper o-edges, the vertices at either side of each "leaf" of the bow tie corresponding to the r-edges, and the vertices at the tips corresponding to the lower o-edges. 
\begin{figure}[H]
	\begin{minipage}{0.48\textwidth}
	\centering
	\includegraphics[width=0.5\textwidth, angle=270]{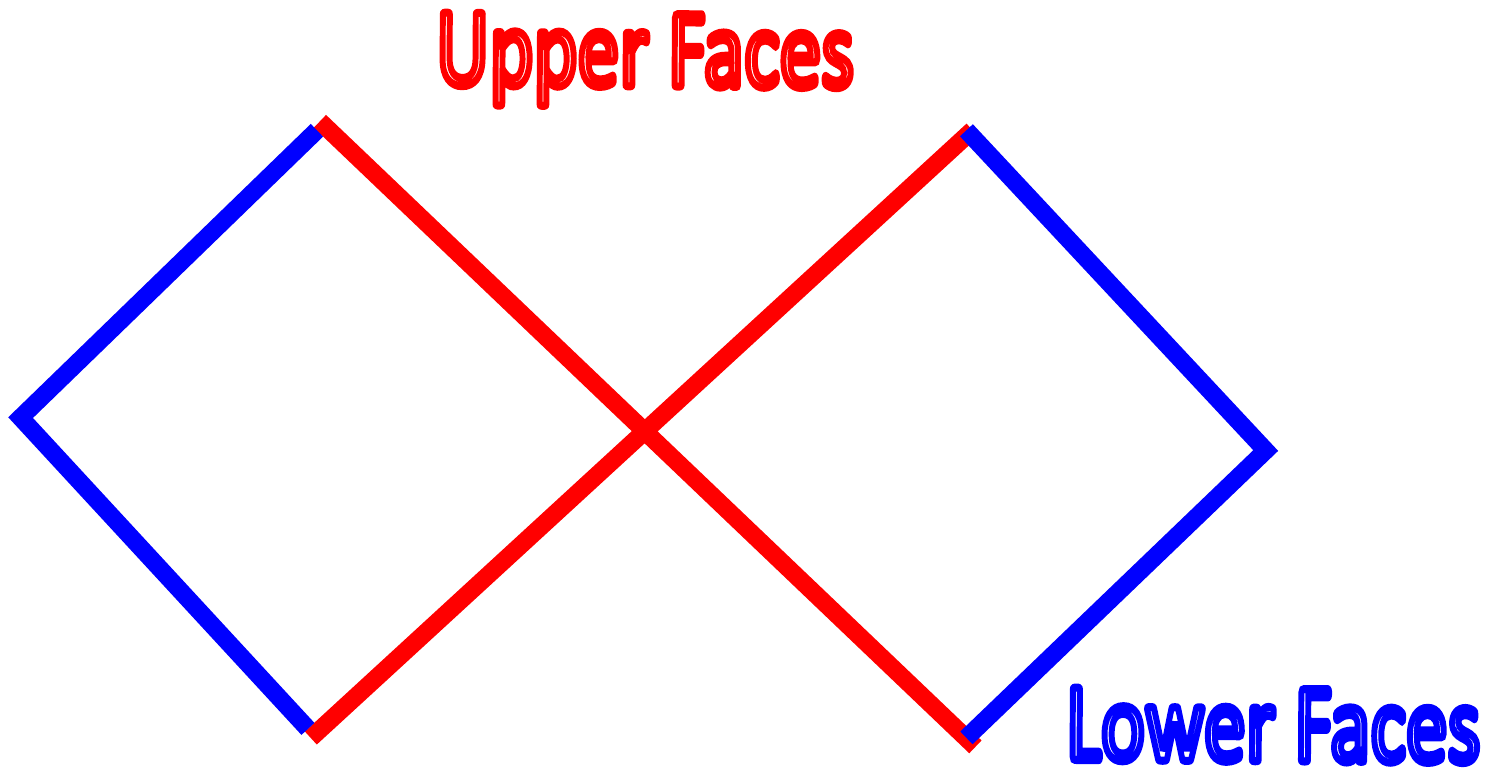}
	\caption{A "Bow Tie" of $\mathcal{V}'(+\infty)$}
	\label{fig:bowtie}
	\end{minipage}\hfill
	\begin{minipage}{0.48\textwidth}
	\centering
	\includegraphics[width=.7\textwidth]{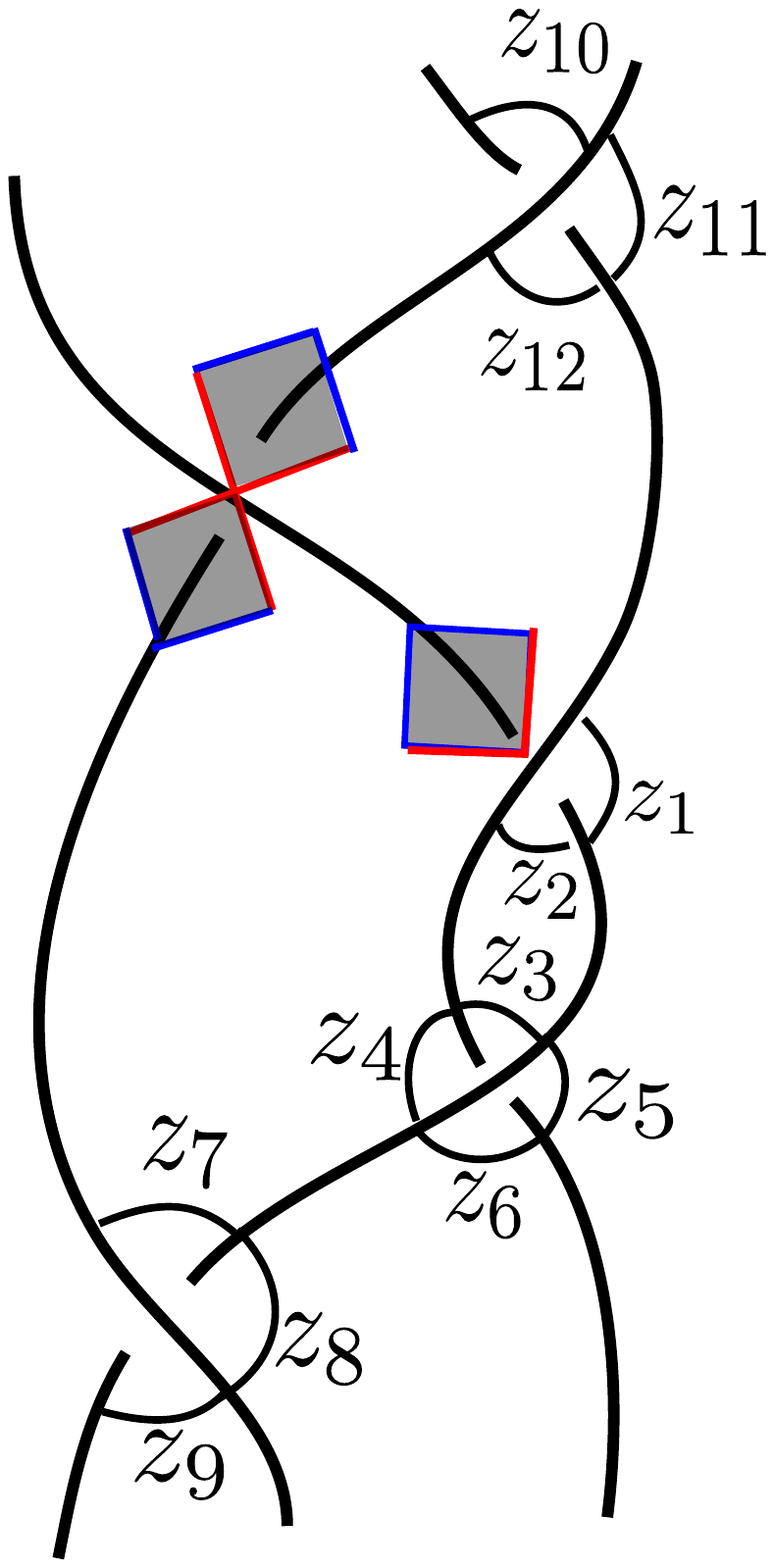}
	\caption{Example}
	\label{fig:topofknot}
	\end{minipage}
\end{figure}
These bowties are glued together, and since our knot is alternating, we will always have inner edges glued to outer edges, where we have all tetrahedra after collapsing. 
\par
Around the $\infty$-segment $\mathcal{V}'(+\infty)$ is missing a bow tie at one crossing, and a half of a bow tie at another crossing, as shown in Figure \ref{fig:topofknot}. Other than that, we have a gluing of top faces to bottom faces around the regions adjacent to the $\infty$-segment, where the vertex in Figure \ref{fig:vposinftri} around which we have our edges glued corresponds to the R-edge of the region.
\begin{figure}[!htb]
	\begin{minipage}{0.48\textwidth}
		\centering
		\includegraphics[width=.6\textwidth]{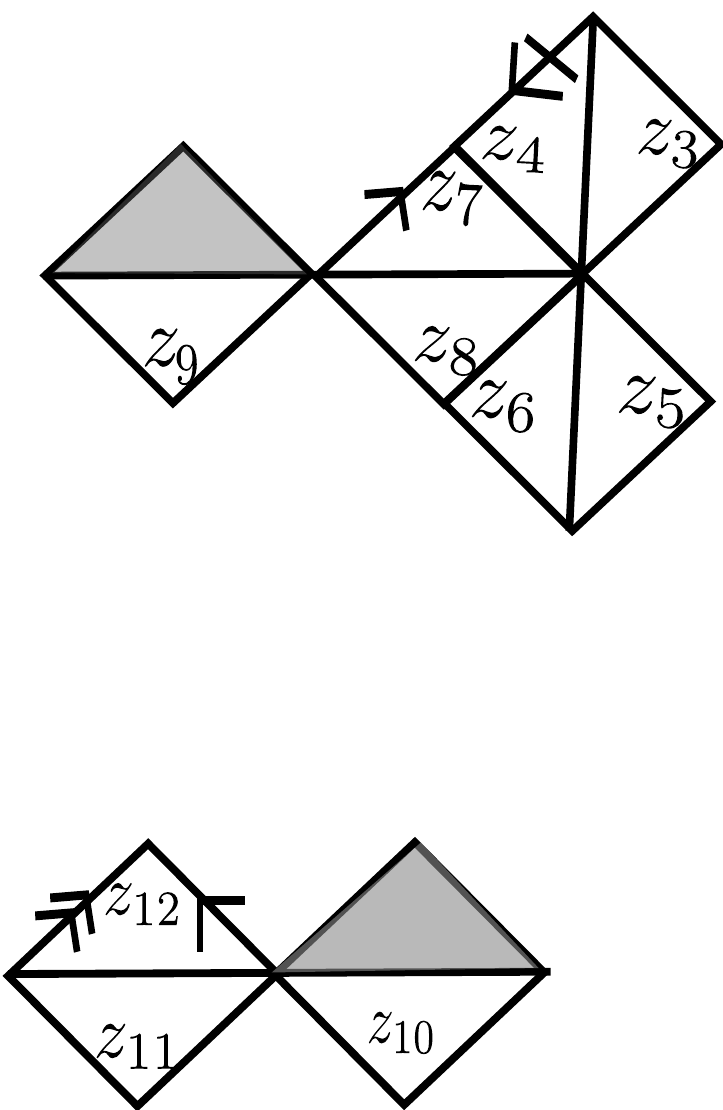}
		\caption{Partial Triangulation of $\mathcal{V}'(+\infty)$}
		\label{fig:vposinftri}
	\end{minipage}\hfill
	\begin{minipage}{0.48\textwidth}
		\centering
		\includegraphics[width=.4\textwidth]{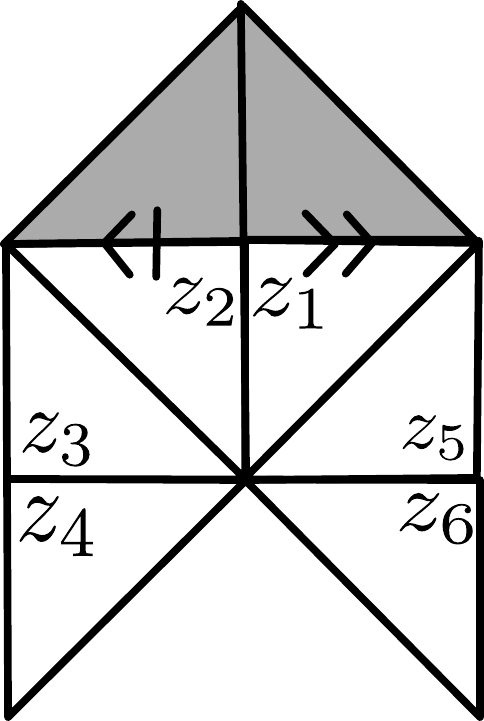}
		\caption{Partial Triangulation of $\mathcal{V}'(K)$}
		\label{fig:vktri}
	\end{minipage}
\end{figure}
Additionally, we have that under the collapsing, two upper faces of a tetrahedron, one in each region, pair to each other, so we have their corresponding edges in $\mathcal{V}'(+\infty)$ glued to each other. The pairing identifies edges of the same type (see \cite{garoufalidis20261} section on collapsing), so the gluing in $\mathcal{V}'(+\infty)$ is such that the center vertices of each bow tie are identified. Finally, we have one lower face in each region paired to an interior face of the octahedron at the crossing where we have two tetrahedra missing. This corresponds to a gluing of edges in $\mathcal{V}'(+\infty)$ with edges in $\mathcal{V}'(K)$, so let us consider the triangulation of $\mathcal{V}'(K)$.
\par
Here we have glued cylinders that when "cut" from the bottom looks like a "flat bow tie" for where an octahedron lies around an over-strand and a "diamond" for where an octahedron lies around an under-strand, and opposite edges are identified. The octahedron with its missing half yields a triangle shape that is half of the diamond, and the flat edges at the top of the diamond shape are the ones that are identified with the edges in $\mathcal{V}'(+\infty)$. Here the lower O-edge of the lower face is identified with the C-edge of the central face (again, see \cite{garoufalidis20261}), so the gluing is reversed compared to the gluing of the other edges around the vertex corresponding to the R-edge.
\par
Therefore, when we complete this gluing, we have a meridian equation of the form:
$$1=z^{-1}_{j_1} z^{-1}_{j_2}\prod z_i ,$$ where $i$ indexes over the tetrahedra in the regions adjacent to the $\infty$-segment and $z_{j_1}, z_{j_2}$ are the shape parameters at the crossing with two missing tetrahedra.
\par
As shown, the completeness equation for this meridian always includes the shape parameter of the tetrahedron at the south corner of the crossing above the region in column 2 adjacent to the $\infty$-segment. Therefore, since we have a gluing equation at that crossing, we can rewrite our meridian completeness equation in terms of the shape parameter at the north corner of that crossing, which is the one not eliminated by the gluing equations at C-edges and R-edges, and eliminate it. In the first case (again, see \cite{garoufalidis20261} and \cite{murakami2018volume}), the meridian equation is already in terms of the shape parameter that was not eliminated. If the chosen variable appears in our equation as $z_j^{-1}$, we can simply take the reciprocal of both sides, and the other side will still be 1, so our equation has the form needed.
\end{proof}

\begin{figure}[H]
	\centering
	\includegraphics[width=.35\linewidth]{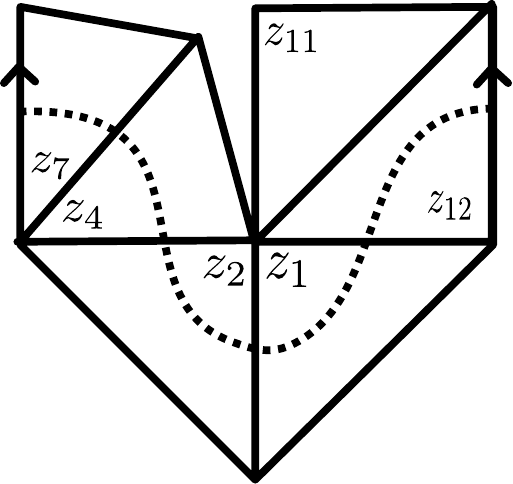}
	\caption{A Meridian of the Knot}
	\label{fig:trimerideq}
\end{figure}

To see an example, consider Figure \ref{fig:topofknot}. Here we have in one region the lower face of $\Delta_7$ paired with the upper face of $\Delta_4$, and their corresponding edges glued in $\mathcal{V}'(+\infty)$. Then, because of our collapsing, we have the upper faces of $\Delta_7$ and $\Delta_{12}$ paired, with their corresponding edges glued. The lower faces of $\Delta_7$ and $\Delta_{12}$ are paired with central faces of $\Delta_2$ and $\Delta_1$, respectively. O-edges of $\Delta_7$ and $\Delta_{12}$ are identified with C-edges of $\Delta_2$ and $\Delta_1$, and R-edges of $\Delta_7$ and $\Delta_{12}$ are identified with lower U-edges of $\Delta_2$ and $\Delta_1$, so the gluing is in accordance with the arrows in the Figures \ref{fig:vposinftri} and \ref{fig:vktri}. Thus we have a completeness equation for the meridian around the dotted line which is $1=z_7z_4z_2^{-1}z_1^{-1}z_{12}$. From the completeness equation at the uppermost crossing, we have $1=z_{10}z_{11}z_{12}$, so we can write this meridian equation as $1=z_{10}(z_7^{-1}z_4^{-1}z_2z_1z_{11}^{-1})$. $z_{10}$ is the shape parameter not eliminated by the gluing equations of the C-edges and R-edges.

\subsection{Determinant Computation}
Now, we know that $$det(G_{\mu}diag(\zeta)+G'_{\mu}diag(\zeta')+G''_{\mu}diag(\zeta''))=det(\frac{\partial (g_1,...,g_{N-1},g_{\mu})}{\partial (z_1,...,z_N)}),$$ where $$g_I=\sum_{i=1}^{N}(G_{Ii}logz_i+G'_{Ii}logz'_i+G''_{Ii}logz''_i),$$ $$g_{\mu}=\sum_{i=1}^{N}(C_i logz_i+C'_i logz'_i+C''_ilogz''_i).$$
\par 
In \cite{dimofte2013quantum}, it is shown that this is independent of choice of which edge equation is replaced by the meridian equation, so it is equal (up to a sign) to $$det(\frac{\partial (g_1,...,g_n,g_{\mu},g_{n+2},...,g_N)}{\partial (z_1,...,z_N)})$$ or $$det(\frac{\partial (g_1,...,g_{n-1},g_{\mu},g_{n+1},...,g_N)}{\partial (z_1,...,z_N)}).$$
We now reorder both the equations and the shape parameters. The shape parameter $z_i$ for $i\in \{1,...,n\}$ will correspond to the corner between segments with labels $x_{i-1}$ and $x_i$ (or $1$ and $x_1$ for $z_1$), and the shape parameters mapped to by $\kappa$ and the meridian equation are $z_{n+1},...,z_N$. The C and R-edge equations are those that appear after $g_{\mu}$, where we choose the first option above if $\kappa$ defines all but one of the eliminations needed, and the second option if $\kappa$ defines all of the eliminations. Now applying Lemma \ref{lem_elim} for all C and R-edge equations, and possibly the meridian equation, we get that $$det(G_{\mu}diag(\zeta)+G'_{\mu}diag(\zeta')+G''_{\mu}diag(\zeta''))=\frac{1}{\prod_{k=n+1}^{N}z_k}det(\frac{\partial (g'_1,...,g'_n)}{\partial (z_1,...,z_n)}),$$ where $g'_i$ is obtained from $g_i$ by eliminating the $N-n$ shape parameters $z_{n+1},...,z_N$.
\par
From \cite{murakami2018volume}, each gluing equation can be obtained by combining two critical point equations $x_i\partial V/\partial x_i$, and thus $$det(\frac{\partial}{\partial x_i}(x_j\frac{\partial V}{\partial x_j}))=det(\frac{\partial (g'_1,...,g'_n)}{\partial (x_1,...,x_n)})=det(\frac{\partial (g'_1,...,g'_n)}{\partial (z_1,...,z_n)})det(\frac{\partial (z_1,...,z_n)}{\partial (x_1,...,x_n)}).$$
\par
Due to our choice of labeling, we additionally have $z_1=x_1^{\pm 1}$ and $z_i=(x_{i-1}/x_i)^{\pm 1}$. Therefore $\frac{\partial (z_1,...,z_n)}{\partial (x_1,...,x_n)}$ is lower triangular, so its determinant is the product of the diagonal entries, $\prod_{i=1}^{n}\frac{\partial z_i}{\partial x_i}$. If $z_1=x_1$, $$\frac{\partial z_1}{\partial x_1}=1=\frac{z_1}{x_1},$$ and if $z_1=x_1^{-1}$, $$\frac{\partial z_1}{\partial x_1}=-\frac{1}{x_1^2}=-\frac{z_1}{x_1}.$$ Likewise for $2 \leq i\leq n$, if $z_i=\frac{x_i}{x_{i-1}}$, $$\frac{\partial z_i}{\partial x_i}=\frac{1}{x_{i-1}}=\frac{z_i}{x_i},$$ and if $z_i=\frac{x_{i-1}}{x_i}$, $$\frac{\partial z_i}{\partial x_i}=-\frac{x_{i-1}}{x_i^2}=-\frac{z_i}{x_i}.$$ Thus $$det(\frac{\partial (z_1,...,z_n)}{\partial (x_1,...,x_n)})=\pm \frac{\prod_{k=1}^{n}z_k}{\prod_{k=1}^{n}x_k}.$$ Combining these results we get
\begin{align}
	det(x_i\frac{\partial}{\partial x_i}(x_j\frac{\partial V}{\partial x_j})) &= \prod_{k=1}^{n}x_kdet(\frac{\partial}{\partial x_i}(x_j\frac{\partial V}{\partial x_j})) \\
	&=\prod_{k=1}^{n}x_kdet(\frac{\partial (g'_1,...,g'_n)}{\partial (z_1,...,z_n)})det(\frac{\partial (z_1,...,z_n)}{\partial (x_1,...,x_n)}) \\
	&=\prod_{k=1}^{n}z_kdet(\frac{\partial (g'_1,...,g'_n)}{\partial (z_1,...,z_n)}) \\
	&=\prod_{k=1}^{N}z_kdet(G_{\mu}diag(\zeta)+G'_{\mu}diag(\zeta')+G''_{\mu}diag(\zeta'')),
\end{align}
proving theorem \ref{mainthm}.

\section*{Declaration of AI Use}
No artificial intelligence (AI) was used in the making of this paper. All mathematical arguments are the author's own.

\bibliography{references}{}
\bibliographystyle{ieeetr}


\end{document}